\documentclass [a4paper,10pt]{article}
\usepackage[utf8]{inputenc}
\usepackage[T1]{fontenc}
\usepackage[english]{babel}
\usepackage{eucal,enumitem}
\usepackage{graphicx}
\usepackage{amsmath,amssymb,amsthm,}
\usepackage{amsfonts}
\usepackage{mathtools}
\usepackage{tikz}
\usepackage{geometry}
\usepackage{cancel}
\usepackage{multicol}
\usepackage{color}
\usepackage{fancyhdr}
\usepackage{afterpage} %pagina bianca
\usepackage{booktabs}
\usepackage{caption}
\usepackage{hyperref}
\usepackage{cleveref}
\usepackage{url}
\usepackage[pdftex]{pict2e}
\usepackage{adjustbox}
\theoremstyle{definition}
\newtheorem{definition}{Definition}[section]
\theoremstyle{plain}
\newtheorem{lemma}[definition]{Lemma}
\theoremstyle{plain}
\newtheorem{theorem}[definition]{Theorem}
\theoremstyle{plain}
\newtheorem{corollary}[definition]{Corollary}
\theoremstyle{plain}
\newtheorem{proposition}[definition]{Proposition}
\theoremstyle{definition}
\newtheorem{remark}[definition]{Remark}
\theoremstyle{definition}
\newtheorem{example}[definition]{Example}
\theoremstyle{remark}

\usepackage{bookmark}
\usepackage{csquotes}
\usepackage[backend=biber,style=numeric,sorting=nyt]{biblatex}
\usepackage{epigraph}
\theoremstyle{plain}
\newtheorem*{customtheoremA}{\textbf{Theorem A}}
\newtheorem*{customtheoremB}{\textbf{Theorem B}}

\renewbibmacro*{in:}{}  % Rimuove la parola "In"

\DeclareFieldFormat[article]{title}{\textit{#1}} 

\DeclareFieldFormat[article]{journaltitle}{#1} 

\DeclareFieldFormat[article]{volume}{\textbf{#1}} 

\title{SOLVABLE GROUPS IN WHICH EVERY REAL ELEMENT HAS PRIME POWER ORDER}

\bigskip

\author{
\bigskip
  \small{BY} \\ 
 \bigskip
  ALESSANDRO GIORGI
  \smallskip
  \thanks{The research was carried out while the author was enrolled at Università degli Studi di Firenze.}\\
  \textit{e-mail: alessandro.giorgi@edu.unifi.it}
}

\date{}

\begin{document}

\maketitle

\begin{abstract}
\noindent
    We study the finite solvable groups $G$ in which every real element has prime power order. We divide our examination into two parts: the case $\textbf{O}_2(G)>1$ and the case $\textbf{O}_2(G)=1$. Specifically we proved that if $\textbf{O}_2(G)>1$ then $G$ is a $\{2,p\}$-group. Finally, by taking into consideration the examples presented in the analysis of the $\textbf{O}_2(G)=1$ case, we deduce some interesting and unexpected results about the connectedness of the real prime graph $\Gamma_{\mathbb{R}}(G)$. In particular, we found that there are groups such that $\Gamma_{\mathbb{R}}(G)$ has respectively 3 and 4 connected components.
\end{abstract}

\section{Introduction}
The prime graph $\Gamma(G)$ of a finite group $G$, also known as the Gruenberg-Kegel graph of $G$, constitutes an important environment to study the "relations" between the elements of $G$ and more generally to analyze the structure of the group $G$. Such graph is defined in this way: the vertices of $\Gamma(G)$ are the prime divisors of $|G|$ and there is an edge between the vertices $p$ and $q$ if $G$ contains an element of order $pq$.

Reality is an interesting and useful notion to consider in finite group theory. An element $x \in G$ is said to be \textbf{real} if it is $G$-conjugate to its inverse $x^{-1}$, i.e. if there exists $g\in G$ such that $x^g=x^{-1}$. We can then define, by just considering the real elements of $G$, the \textbf{real prime graph} $\Gamma_{\mathbb{R}}(G)$ analogously to the prime graph: the vertices of $\Gamma_{\mathbb{R}}(G)$ are the primes $p$ such that $G$ contains a real element of order $p$ and the vertices $p$ and $q$ are connected if $G$ contains a real element of order $pq$.

In this paper we study the finite solvable groups $G$ such that all the vertices of $\Gamma_{\mathbb{R}}(G)$ are isolated, starting from the results obtained by Dolfi, Gluck and Navarro in \cite{DGN}, where the authors investigated the finite solvable groups for which $2$ is an isolated vertex of $\Gamma_{\mathbb{R}}(G)$. Taking into consideration the above definition of $\Gamma_{\mathbb{R}}(G)$, it is clear that the fact that the vertices of $\Gamma_{\mathbb{R}}(G)$ are all isolated is equivalent to the fact that every real element of $G$ has prime power order. The analogous problem for $\Gamma(G)$ was studied in 1957 by G. Higman in \cite{Hi}, in which he considered both solvable and insolvable groups $G$ for which every element has prime power order. For the solvable ones, his main result contained in [Theorem 1, \cite{Hi}] tells us that $G$ must be a $p$-group (that could be considered as the trivial case for this question) or a $\{p,q\}$-group ([Theorem 1, \cite{Hi}] actually says a lot more about the structure of $G$).

It is easy to see that if $G$ contains a real element other than the identity then $|G|$ is even. Also, since every involution is a real element, we can conclude that if $G$ has non-trivial real elements then $2\mid |G|$ and $2$ is a vertex of $\Gamma_{\mathbb{R}}(G)$. So what we initially wanted to prove, similarly to Higman, was that if every real element in $G$ has prime power order then $G$ must be either a $2$-group or a $\{2,p\}$-group. However, as we will see, that is not true in general. Still we were able to prove it in the case that $\textbf{O}_2(G)>1$, where $\textbf{O}_2(G)$ is the largest normal $2$-subgroup of $G$.

On a more technical side, we note that since the real elements of $G$ are the real elements of $\textbf{O}^{2'}(G)$, the smallest normal subgroup of $G$ with odd index, when investigating the real elements of a finite group we can assume that $G=\textbf{O}^{2'}(G)$.

With that being said, our first main result is the following:

\begin{customtheoremA}
Suppose that $G$ is a finite solvable group with $\textbf{O}^{2'}(G)=G$. Let $N=\textbf{O}_2(G)$. Suppose also that $G$ is not a $2$-group and $N>1$. Then the following are equivalent:
\begin{enumerate}
    \item every real element of $G$ has prime power order;
    \item $G$ is a $\{2,p\}$-group, with $p$ odd prime, $G=N\rtimes (K\rtimes Q)$ is a $2$-Frobenius group, with $K$ a cyclic $p$-group and $Q$ a cyclic $2$-group, and every element of $G$ has prime power order. In particular if $|Q|=2$ then $K\rtimes Q$ is a dihedral group. In any case, every element of $K$ is real and inverted by $z$, where $z$ is the only involution of $Q$. 
\end{enumerate}
\end{customtheoremA}

But, as we already mentioned above, it is not generally true that if every real element of $G$ has prime power order then $G$ must be a $2$-group or a $\{2,p\}$-group. As a matter of fact, studying the case $\textbf{O}_2(G)=1$, even though we were not able in this paper to precisely describe the structure of such groups, we found examples of groups $G$ such that every real element has prime power order and $|G|$ has respectively three and four prime divisors. We should now recall that by [Corollary B, \cite{DMN}], since we are assuming $G$ solvable, for every prime $p$ such that $p\mid |G|$ we have that $G$ contains a real element of order $p$, so $p$ is a vertex of $\Gamma_{\mathbb{R}}(G)$ by definition. 
In the framework of the real prime graph $\Gamma_{\mathbb{R}}(G)$, since every isolated vertex obviously is a connected component, if we write $n(\Gamma)$ to indicate the number of connected components of a graph $\Gamma$, then we can state our second main result in the following way:

\begin{customtheoremB}
 There exist finite solvable groups $G$ and $H$ with $n(\Gamma_{\mathbb{R}}(G))=3$ and $n(\Gamma_{\mathbb{R}}(H))=4$.

\end{customtheoremB}

The content of Theorem B is somewhat unexpected because it goes against the usual similarity of properties that $\Gamma_{\mathbb{R}}(G)$ has in respect to other two notorious graphs associated to $G$, the prime graph on real character degrees $\Gamma_{cd,\mathbb{R}}(G)$ and the prime graph on real class sizes $\Gamma_{cs,\mathbb{R}}(G)$, for which it is known that $n(\Gamma_{cd,\mathbb{R}}(G))\leq 2$ and $n(\Gamma_{cs,\mathbb{R}}(G))\leq 2$, as we will explain better later. 

It could also be an interesting topic of study to investigate how much the number of prime divisors of $|G|$ can be increased while preserving the condition that every real element of $G$ has prime power order, which is closely related, even though not equivalent, to the open problem of determining the least upper bound for $n(\Gamma_{\mathbb{R}}(G))$.

\medskip

%The author would like to thank Professor S. Dolfi for bringing this problem to his attention and for the invaluable advice received while working under his supervision on this research.

\section{Preliminary results}
For our study we will need some basic but nevertheless fundamental properties of real elements. Let us start with the following lemma (see Lemma 3.2 of \cite{DMN} for a proof).

\begin{lemma}\label{lem:elem prop}
    Let $G$ be a finite group.
    \begin{enumerate}
        \item If $x\in G$ is real, then there is a $2$-element $y\in G$ such that $x^y=x^{-1}$.
        \item If $x \in G$ is real, then $x^m$ is real for every integer $m$.
        \item Suppose that $N \trianglelefteq G$ and that $xN\in G/N$ is real. Suppose also that $o(xN)$ is odd. Then there is a real $y \in G$ such that $xN=yN$.
        \item If $Q$ is a $2$-group acting non-trivially on $G$, then there are $1\neq x\in G$ and $q\in Q$ such that $x^q=x^{-1}$.
    \end{enumerate}
\end{lemma}

The next lemma will be used in several coming proofs to conclude various reasoning by contradiction, since it gives us a sufficient condition (actually it is also a necessary condition, but we will use it only in one "direction") for the existence of real elements of non prime power order. For this reason, even if it is fairly easy, we shall prove it.

\begin{lemma}\label{lem: suff cond}
Let $G$ be a finite group and let $x,y \in G$ be real elements such that $o(x)=p$ and $o(y)=q$, with $p,q$ primes and $p\neq q$. Suppose that $x$ and $y$ are inverted by the same $g\in G$ and that $xy=yx$. Then $xy$ is real inverted by $g$ and $o(xy)=pq$.
\begin{proof}
     We start by proving that $xy$ is real inverted by $g$:
    \begin{equation*}
        (xy)^g=x^gy^g=x^{-1}y^{-1}=y^{-1}x^{-1}=(xy)^{-1},
    \end{equation*}
    since if $x$ and $y$ commute, then also do their inverses. 
    
    Moreover $(xy)^{pq}=1$ which implies $o(xy)\mid pq$. It is clear that $o(xy)\neq 1,p,q$ and so we have $o(xy)=pq$. 
\end{proof}
    
\end{lemma}

We note that \Cref{lem: suff cond} is still valid even if $o(x)=p^{\alpha}$ and $o(y)=q^{\beta}$, but it is sufficient to enunciate it in this form, since by (2) of \Cref{lem:elem prop} we can always consider appropriate powers of $x$ and $y$ and get real elements with prime order.

In continuity with \cite{DGN}, we say that a finite group $G$ \textbf{satisfies R} if every real element has $2$-power order or $2'$-order, and we say that $G$ \textbf{satisfies P} if every real element has prime power order.

We want now to prove that our working hypotheses descend to the quotient. It is well known that solvability does that and it is not difficult to verify that if $G$ is a finite group with $\textbf{O}^{2'}(G)=G$ then, if $N\trianglelefteq G$, we have $\textbf{O}^{2'}(G/N)=G/N$. The property $\textbf{P}$ requires a little bit more work. 

\begin{lemma}\label{lem: P quotient}
Let $G$ be a finite group and let $N\trianglelefteq G$. If $G$ satisfies \textbf{P}, then $G/N$ satisfies \textbf{P}.

\begin{proof}
Since $G$ satisfies \textbf{P} then obviously $G$ satisfies \textbf{R}. So by Lemma 2.2 of \cite{DGN} we have that $G/N$ satisfies \textbf{R}, that is every real element of $G/N$ has $2$-power order or odd order. Suppose by contradiction that $G/N$ contains a real element $xN$ of odd non prime power order. Then by $(3)$ of \Cref{lem:elem prop} there exist a real $y\in G$ such that $xN=yN$. Then we have $o(xN)\mid o(y)$ and so $y$ is a real element of $G$ with non prime power order, against assumptions.
    
\end{proof}
    
\end{lemma}

It is worth noting that \Cref{lem: P quotient} could also be proved directly by extending Lemma 2.1 of \cite{DGN} to our case, instead of integrating Lemma 2.2 of \cite{DGN} with (3) of \Cref{lem:elem prop} as we did above.

We are now going to enunciate Theorem A of \cite{DGN}, which establishes the basis for our work.

\begin{theorem}\label{theo: main DGN}
     Suppose that $G$ is a finite solvable group with $\textbf{O}^{2'}(G)=G$. Assume that every real element of $G$ is either a $2$-element or a $2'$-element. Let $N=\textbf{O}_2(G)$ and $Q\in Syl_2(G)$, and assume that $G$ is not a $2$-group. Then:
    \begin{enumerate}
        \item $G/N$ has a normal $2$-complement $K/N$ and $Q/N$ is cyclic or generalized quaternion. If $zN$ is the unique involution of $Q/N$, then 
        \begin{equation*}
         C_{K/N}(Q/N)=C_{K/N}(zN).
        \end{equation*}
        \item Suppose that $N>1$. Then $N=\textbf{F}(G)$, $Q/N$ is cyclic and $G$ splits over $N$. If $|Q/N|>2$, then $K/N$ is cyclic. In any case, $K/F_2$ is metabelian and $F_2/N$ is abelian, where $F_2/N=\textbf{F}(G/N)$. If $|G|$ is coprime to $3$, then $K/F_2$ is abelian.
    \end{enumerate}
\end{theorem}

We are now going to see some lemmas that give us some initial information on the structure of the group $G/N$ described in \Cref{theo: main DGN}. For the sake of brevity, we are going to write $G$ instead of $G/N$, $K$ instead of $K/N$ and so on.

\begin{lemma}\label{lem: interder KQ}
    Let $G$ be a finite group with $G=K\rtimes Q$ and $\textbf{O}^{2'}(G)=G$, where $K$ is the normal $2$-complement and $Q\in Syl_2(G)$. Then we have $[K,Q]=K$.
    \begin{proof}
        Since $K\trianglelefteq G$, then $[K,Q]\leq K$ and we also know that $[K,Q]\trianglelefteq \langle K,Q \rangle=G$. We can consider $L=[K,Q]Q$. We have that $L\trianglelefteq G$, since $[K,Q]\trianglelefteq G$ and $Q^g=Q^{qk}=Q^k\subseteq [K,Q]Q$ for every $g=qk \in G$, with $q\in Q$ and $k\in K$. Then $L$ is a normal subgroup of $G$ with odd index and so it must be $L=G$. Then it follows that $[K,Q]=K$.
    \end{proof}
    
\end{lemma}

\begin{lemma}\label{lem: abel K}
Let $G$ be a finite group with $G=K\rtimes Q$, where $Q\in Syl_2(G)$ and $K$ is the normal 2-complement of $G$. Suppose that $\textbf{O}^{2'}(G)=G$ and that $Q$ has a unique involution $z$, and assume that $C_K(z)=C_K(Q)$. Then $K$ is abelian if and only if $C_K(z)=1$ (if and only if $z$ inverts every element of $K$).
\begin{proof}
 Suppose $K$ is abelian. Since $(|Q|,|K|)=1$, by coprime action of $Q$ on $K$ we have $K=[K,Q]\times C_K(Q)$. But by \Cref{lem: interder KQ} we know that $[K,Q]=K$ and so $C_K(Q)=1$. It follows that $C_K(z)=C_K(Q)=1$.

 Conversely suppose that $C_K(z)=1$. So $z$ induces an automorphism on $K$ of order $2$ with no fixed points. It is a well known result that then $z$ acts on $K$ as the inversion. Since it is easy to see that a group which has the inversion as an automorphism is abelian, we conclude.

\end{proof}
    
\end{lemma}

Let us finish this section with a remark that contains an idea which will be used several times later on.
\begin{remark}\label{rem: nil impl p group}

Let $G$ be a finite group as in $(1)$ of \Cref{theo: main DGN}, or as in \Cref{lem: abel K}, and assume also that $G$ satisfies \textbf{P}. We prove that if the normal $2$-complement $K$ is nilpotent then $K$ is a $p$-group, for some odd prime $p$.

Suppose $K$ is nilpotent. Assume by contradiction that there are at least two odd primes $p,q$, with $p\neq q$, such that $p,q\mid |K|$. By nilpotence, we have that the $p$-elements commute with the $q$-elements. Consider now the actions of $z$, which is as usual the involution of $Q$, on $\textbf{O}_p(K)$ and on $\textbf{O}_q(K)$. If $z$ centralizes, to fix ideas, $\textbf{O}_p(K)$, then $\textbf{O}_p(K)\leq C_K(z)=C_K(Q)$. So we would have $Q\leq C_G(\textbf{O}_p(K))$. %Then $C_G(\textbf{O}_p(K))$ would be a normal subgroup of $G$ with odd index, which would imply $C_G(\textbf{O}_p(K))=G$, by the hypothesis $\textbf{O}^{2'}(G)=G$. Then $\textbf{O}_p(K)\leq \textbf{Z}(G)$. 
In particular, $\textbf{O}_{p'}(K)\rtimes Q$ would be a normal subgroup of $G$ with odd index, against the fact that $\textbf{O}^{2'}(G)=G$. The same goes for $\textbf{O}_q(K)$. So the actions of $z$ on both are not trivial. Then by $(4)$ of \Cref{lem:elem prop} there exist real elements $x\in \textbf{O}_p(K)$ and $y\in \textbf{O}_q(K)$ inverted by $z$, respectively of order $p$ and $q$, that commute. But then by \Cref{lem: suff cond} $xy$ is a real element, inverted by $z$, of order $pq$, against \textbf{P}. Hence $K$ must be a $p$-group.

\end{remark}

\section{The case $\textbf{O}_2(G)>1$}
In order to treat this case it is necessary to repeat some ideas and results introduced in \cite{DGN}. 

\begin{definition} (Standard Hypotheses)
Let $G=KQ$, where $K>1$ is normal of odd order, $Q\in Syl_2(G)$ is cyclic or generalized quaternion and $C_K(Q)=C_K(z)$, with $z$ the unique involution of $Q$. Suppose also that $\textbf{O}^{2'}(G)=G$. Assume that $G$ acts on a $2$-group $V$ and that $C_G(v)$ has a normal Sylow $2$-subgroup for all $1\neq v\in V$. In this case, we say that $G$ satisfies the Standard Hypotheses with respect to $V$.
\end{definition}

The following theorem (see [Theorem 3.1, \cite{DGN}]) explains the introduction of such hypotheses.

\begin{theorem}\label{theo: SH}
Suppose $G$ is a finite solvable group with $\textbf{O}^{2'}(G)=G$. Assume that $G$ satisfies \textbf{R} and that $G>N=\textbf{O}_2(G)>1$. Then there exists a subgroup $H$ of $G$ such that $G=NH$, with $N\cap H=1$, such that $H$ satisfies the Standard Hypotheses with respect to $N$. Moreover $\textbf{O}_{2'}(G)\leq \textbf{Z}(G)$. 

Conversely, if $H=KQ$ satisfies the Standard Hypotheses with respect to $V$, then $G=V\rtimes H$ satisfies \textbf{R}.
    
\end{theorem}

We would like now to have a set of assumptions that descend to quotients. Recall that, in a group action of $Q$ on $K$, $A/B$ is said a $Q$-\textit{invariant} $p$-\textit{section}, for a prime $p$, if $A,B$ are $Q$-invariant subgroups of $K$, $B\trianglelefteq A$ and $A/B$ is a $p$-group.

\begin{definition} (H2)
Let $G$ be a group with a cyclic Sylow $2$-subgroup $Q>1$ and a normal $2$-complement $K$. Suppose that $\textbf{O}^{2'}(G)=G$ and $C_K(Q)=C_K(z)$, where $z$ is the unique involution of $Q$. Assume also that, for every prime $p$ and for every $Q$-invariant $p$-section $A/B$ of $K$, $[A/B,Q]$ is cyclic.
    
\end{definition}

Then we have that the hypotheses (H2) descend to quotients on $N\trianglelefteq G$, $N\leq K$ (see [Lemma 3.2, \cite{DGN}]).

\begin{lemma}\label{lem: H2 quotient}
  Suppose that $G$ satisfies \textnormal{(H2)} and let $N\trianglelefteq G$, $N\leq K$. Then also $G/N$ satisfies \textnormal{(H2)}.  
\end{lemma}

The following proposition clarifies the connection between the Standard Hypotheses and (H2) (see [Proposition 3.4, \cite{DGN}]).
\begin{proposition} \label{SH implies H2}
   If $G=KQ$ satisfies the Standard Hypotheses, then $G$ satisfies \textnormal{(H2)}. 
\end{proposition}

Before we can start with our own investigation, we need other two results from \cite{DGN}, which describe the chief factors and the structure of the Sylow subgroups of a group $G$ satisfying (H2) (see respectively [Proposition 3.7, \cite{DGN}] and [Theorem 3.9, \cite{DGN}]).

\begin{proposition}\label{prop: chief fact}
Let $G$ satisfy \textnormal{(H2)}. Let $p$ be a prime such that $p\mid |K|$ and let $P=\textbf{O}_p(G)$, $R=\Phi(G)\cap P$ and $X=P/R$. If $X\neq 1$, then $X$ is a noncentral chief factor of $G$, $P\in Syl_p(G)$ and $R=\Phi(P)$.
    
\end{proposition}

\begin{theorem}\label{theo: structure Syl}
Let $G$ satisfy \textnormal{(H2)} with $K>1$. Let $p$ be a prime such that $p\mid |K|$ and $P\in Syl_p(G)$. Then $P$ is homocyclic abelian of rank at most $3$. If $p$ divides $|K/K'|$, then $P$ is cyclic. Also we have $\textbf{Z}(G)=1$.
    
\end{theorem}

The first consequence that we deduce is the following:
\begin{proposition}\label{prop: Fitting is p group}
Let $G$ satisfy \textnormal{(H2)} and assume $G$ satisfies \textbf{P}, with $K>1$. Then $\textbf{F}(G)=\textbf{F}(K)$ and $\textbf{F}(G)\in Syl_p(G)$, for $p$ (odd) prime such that $p\mid |K|$.

\begin{proof}
    Let us start by proving that $\textbf{F}(G)=\textbf{F}(K)$. To do so, it is sufficient to show that $\textbf{F}(G)$ has odd order.  Indeed in that case $\textbf{F}(G)\leq K$ and so $\textbf{F}(G)\leq \textbf{F}(K)$. The other incusion is trivial. Suppose then by contradiction that $\textbf{F}(G)$ has even order. Then it must be $\textbf{O}_2(G)>1$ and so $z \in \textbf{O}_2(G)$, where $z$ is the unique involution of $Q$. Therefore $z$ commutes with $\textbf{O}_{2'}(G)=K$, that is $K=C_K(z)=C_K(Q)$. Then $Q\trianglelefteq G$ and $G=K\times Q$. Since $\textbf{O}^{2'}(G)=G$, we have $G=Q$ and $K=1$, against the assumptions. 

    Working by contradiction, assume now that $\textbf{F}(G)$ is not a $p$-group. Then there exist primes $p,q$, $p\neq q$ such that $\textbf{O}_p(G),\textbf{O}_q(G)>1$. Let us consider the action of $z$ on $\textbf{O}_j(G)$, for $j=p,q$. If the action is trivial, then, as in \Cref{rem: nil impl p group}, we have $Q\leq C_G(\textbf{O}_j(G))$. So $C_G(\textbf{O}_j(G))$ is a normal subgroup with odd index. Therefore, since $\textbf{O}^{2'}(G)=G$, we have that $C_G(\textbf{O}_j(G))=G$, that is $\textbf{O}_j(G)\leq \textbf{Z}(G)$, against the fact that $\textbf{Z}(G)=1$ by \Cref{theo: structure Syl}. So both those actions are not trivial. Thus by (4) of \Cref{lem:elem prop} there exist a real $p$-element $x$ and a real $q$-element $y$ both inverted by $z$, and they commute. Then by \Cref{lem: suff cond} $xy$ is a real element of order $pq$, against the property \textbf{P}. So $\textbf{F}(G)=\textbf{O}_p(G)$, for some (odd) prime $p$.

    Finally, we note that a group satisfying (H2) is solvable by the Feit-Thompson theorem. It is well known that in a finite non-trivial solvable group $\textbf{F}(G)>1$ and $\textbf{F}(G)>\Phi(G)$. Then, since $R=\textbf{F}(G)\cap \Phi(G)=\Phi(G)$, we have $\textbf{F}(G)/R=\textbf{F}(G)/\Phi(G)>1$. Therefore, by \Cref{prop: chief fact} we conclude that $\textbf{F}(G)\in Syl_p(G)$. 

\end{proof}
    
\end{proposition}

We are almost ready to prove the main theorem of this paper. We just need to simplify our work with one last observation. Suppose $G$ is a finite solvable group, with $\textbf{O}^{2'}(G)=G$, $G>N=\textbf{O}_2(G)>1$ and assume that $G$ satisfies \textbf{P}. In particular, we can apply \Cref{theo: main DGN} to $G$. Since we are assuming $N>1$, if we also suppose $|Q/N|>2$, then by (2) of \Cref{theo: main DGN} we have that $K/N$ is cyclic. Therefore, by \Cref{rem: nil impl p group}, or even with little more work by \Cref{lem: abel K}, we know that $K/N$ is a $p$-group.

So it suffices to consider the case $|Q/N|=2$, for which we have:

\begin{theorem}\label{theo: case |Q|=2}
Let $G$ be a finite solvable group, with $\textbf{O}^{2'}(G)=G$, and suppose that $G$ satisfies \textbf{P}. Let $N=\textbf{O}_2(G)$ and $Q\in Syl_2(G)$. Assume that $G$ is not a $2$-group, $N>1$ and $|Q/N|=2$. Then we have $G=N\rtimes H$, $H=K\rtimes \langle z \rangle$, with $\langle z \rangle \in Syl_2(H)$ and $o(z)=2$, and $K$ is a cyclic $p$-group, for some (odd) prime $p$. Moreover $H$ is a dihedral group. In particular, every element of $K$ is real in $H$, inverted by $z$.

\begin{proof}
    Since every hypothesis descend to quotients, we may assume that $N$ is minimal normal in $G$, so that it is an elementary abelian $2$-group. Therefore all his elements are involutions, and real elements of $G$, and $[N,N]=1$.

    By \Cref{theo: SH} we know that there exists $H\leq G$ such that $G=N\rtimes H$ and $H$ satisfies the Standard Hypotheses with respect to $N$. So $H=K\rtimes Q_0$, with $K>1$ the normal $2$-complement of $H$ and $Q_0\in Syl_2(H)$. But $H\cong G/N$, therefore $|Q_0|=2$. So $H=K\rtimes \langle z \rangle$, with $z$ involution. Furthermore, by \Cref{SH implies H2} we have that $H$ satisfies (H2).

    Working by contradiction, assume that $K$ is not a $p$-group. Then $K$ is not nilpotent, since otherwise, as in \Cref{rem: nil impl p group}, $K$ would be a $p$-group. Let us consider $K^{\infty}$ the residual nilpotent of $K$, which for a finite group is the last term of the lower central series and the smallest normal subgroup such that the quotient is nilpotent. We have $K^{\infty} \,char \,K$ and so $K^{\infty} \trianglelefteq H$. Moreover $K^{\infty}\neq 1$. Thus $H/K^{\infty}$ has the normal $2$-complement $K/K^{\infty}$, which is nilpotent, and so $K/K^{\infty}$ is a $p$-group. Also, by \Cref{lem: interder KQ}, we have $[K/K^{\infty},\langle z \rangle K^{\infty}/K^{\infty}]=K/K^{\infty}$. Then $K/K^{\infty}=[K/K^{\infty},z]$, which is cyclic since $H$ satisfies (H2). In particular, $K/K^{\infty}$ is abelian, therefore $K'\leq K^{\infty}$ and we conclude $K'=K^{\infty}$.

    At this point, by the results in (2) of \Cref{theo: main DGN}, it is natural to split the proof into two parts: the case $3\nmid
    |K|$, which will be part (a), and the case $3\mid |K|$, which will be part (b).

\medskip
(a) Assume $3\nmid |K|$. Then by \Cref{theo: main DGN} we know $K/\textbf{F}(H)$ is abelian. But $\textbf{F}(H)=\textbf{F}(K)$, since $H\cong G/N$ and so $\textbf{O}_2(H)=1$. Thus $K/\textbf{F}(K)$ is abelian and therefore $K'\leq \textbf{F}(K)$. Furthermore, by \Cref{prop: Fitting is p group} $\textbf{F}(K)$ is a $q$-group, for some prime $q$ with $q\mid |K|$. Since we are assuming that $K$ is not a $p$-group, it must be $p\neq q$. Then also it must be $K'=\textbf{F}(K)$, otherwise we would have $p\mid |\textbf{F}(K)|$. For convenience, let $L=K^{\infty}=K'=\textbf{F}(K)$. We note that $L\in Syl_q(H)$, so by \Cref{theo: structure Syl} $L$ is homocyclic of rank at most $3$. From what has been said so far, we have $K=L\rtimes P$, where $P\in Syl_p(K)$, $P\cong K/L$ and so $P$ is cyclic. Moreover, up to conjugacy, we can assume $z\in N_H(P)$. Indeed, by the Frattini argument we have $H=KN_H(P)$ and so 
 \begin{equation*}
            2|K|=|H|=|KN_H(P)|=\frac{|K||N_H(P)|}{|K\cap N_H(P)|}=\frac{|K||N_H(P)|}{|N_K(P)|},
        \end{equation*}
so that $|N_H(P)|=2|N_K(P)|$. Then we can assume that $z$ acts on $P$. We see that $C_P(z)=1$, otherwise in $H/L\cong P \langle z \rangle$ we would have $K/L\neq [K/L,z]$ since $C_{K/L}(z)\neq 1$ (recall $L=K^{\infty}$). Therefore $z$ acts fixed-point free on $P$ and it is well known that then $z$ invertes all the elements of $P$. So $P\langle z \rangle$ is dihedral and in particular it is a Frobenius group.

We want now to verify that $C_L(P)=1$, which follows easily by proving $\textbf{Z}(K)=1$. Assume by contradiction $1\neq \textbf{Z}(K)$. Obviously $\textbf{Z}(K)\trianglelefteq H$. Consider then the action of $z$ on $\textbf{Z}(K)$. Such action is not trivial, otherwise we would have, since $\textbf{O}^{2'}(H)=H$, $1\neq \textbf{Z}(K)\leq \textbf{Z}(H)$, against the fact that $\textbf{Z}(H)=1$ by \Cref{theo: structure Syl}. So, by (4) of \Cref{lem:elem prop}, there is a real element $1\neq x\in \textbf{Z}(K)$, $x=ly, l\in L, y\in P$, inverted by $z$. If $l\neq 1$, then an approriate power of $x$ would be a real $q$-element inverted by $z$ that commutes with every element of $P$. Recalling that every element of $P$ is a real $p$-element inverted by $z$, by \Cref{lem: suff cond} we would get a real element of non prime power order, against \textbf{P}. So $l=1$ and $x$ is a real $p$-element inverted by $z$ that commutes with every element of $L$. But $z$ acts non-trivially on $L$ (otherwise $1<L\leq \textbf{Z}(H)$) and so $L$ contains at least one non-trivial real $q$-element inverted by $z$. As before, that goes against \textbf{P}. 

We can now prove that $L$ has precisely rank $2$. We have $\Phi(L)\trianglelefteq H$ (actually by \Cref{prop: chief fact} $\Phi(L)=\Phi(H)$). So $P\langle z \rangle$ acts on $\Phi(L)$. Let us consider the quotient $L/\Phi(L)$. Using the bar notation, we have that $\overline{L}$ is a $\mathbb{Z}_q$-vector space of dimension at most 3, by \Cref{theo: structure Syl}. By coprime action, we have $C_{\overline{L}}(P)=\overline{C_L(P)}=1$. Recall also that $P\langle z \rangle$ is a Frobenius group, with kernel $P$ and complement $\langle z \rangle$. So, considering $\overline{L}$ as a $\mathbb{Z}_q[P\langle z \rangle]$-module, by [Theorem 15.16, \cite{I}] we deduce $\dim_{\mathbb{Z}_q}(\overline{L})=|\langle z \rangle| \dim_{\mathbb{Z}_q}(C_{\overline{L}}(\langle z \rangle))$. In particular $2\mid \dim_{\mathbb{Z}_q}(\overline{L})$ and $\dim_{\mathbb{Z}_q}(\overline{L})\leq 3$, so it must  be $\dim_{\mathbb{Z}_q}(\overline{L})=2$. Then by Burnside basis theorem $L$ has rank 2.

\vspace{1.5cm}
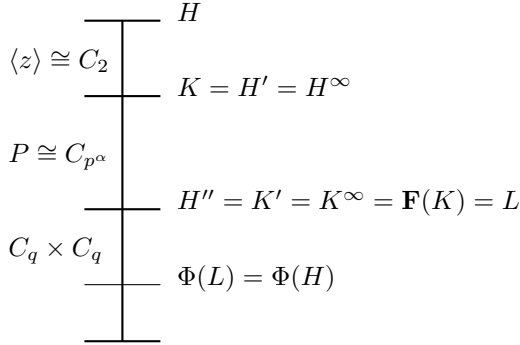
\begin{figure}[ht]
\setlength{\unitlength}{0,5cm}
\begin{picture}(20,20)(-12.5,-12.5)
\put(-1,3){\line(1,0){2}}
\thicklines
\put(0,1.5){\line(0,1){8.5}}
\put(-1,10){\line(1,0){2}}
\put(1,10){\hbox{\kern3pt\ $H$}}
\put(-1,8){\line(1,0){2}}
\put(1,8){\hbox{\kern3pt\ $K=H'=H^{\infty}$}}
\put(-1,5){\line(1,0){2}}
\put(1,5){\hbox{\kern3pt\ $H''=K'=K^{\infty}=\textbf{F}(K)=L$}}
\put(1,3){\hbox{\kern3pt\ $\Phi(L)=\Phi(H)$}}
\put (-1,1.5){\line(1,0){2}}
\put(-3,8.75){$\langle z \rangle \cong C_2$}
\put(-3,6.25){$P\cong C_{p^{\alpha}}$}
\put(-3,3.75){$C_q\times C_q$}

\end{picture}
\vspace{-6.8cm}
\caption{Structure of group $H$ in part (a)}
\end{figure}
%\vspace{0.5cm}

At this point, we want to show that $K=L\rtimes P$ is a Frobenius group. Considering the quotient $L/\Phi(L)$ we can assume that $\Phi(L)=1$. By \Cref{prop: chief fact} we then have that $L$ is minimal normal in $H$. Also, if we write $\mathbb{F}=\mathbb{Z}_q$, we have seen above that $L$ is an $\mathbb{F}$-vector space of dimension $2$. So $L$ is an irreducible $\mathbb{F}[H]$-module. Since $K\trianglelefteq H$, by Clifford theorem we know that $L$ is completely reducible as an $\mathbb{F}[K]$-module (which we write as $L_K$). If $L_K$ is irreducible, then it is a faithful and irreducible $\mathbb{F}[P]$-module, and $P$ is cyclic. It is well known that then $P$ acts fixed-point free on $L$. Suppose therefore that $L_K=L_1\oplus L_2$, with $L_1,L_2$ irreducible $\mathbb{F}[K]$-modules. We want to verify that the irreducible components $L_1,L_2$ have the same kernel. If $L_K$ is homogeneous then it is obvious, since $L_1\cong L_2$ as $\mathbb{F}[K]$-modules. Assume that $L_K$ is not homogeneous. By Clifford theorem we know that $G/K\cong \langle z \rangle$ acts transitively on the set of homogeneous components, in this case $\{L_1,L_2\}$. Thus $L_1^{z}=L_2$ and therefore $C_K(L_1)^z=C_K(L_1^z)=C_K(L_2)$, that is the kernels are $z$-conjugate. But $C_K(L_1)=LX$, for some $X\leq P$, and so $C_K(L_2)=C_K(L_1)^z=L^zX^z=LX$, since $L\trianglelefteq H$ and $z$ normalizes every subgroup of $P$. In any case, we proved $C_K(L_1)=C_K(L_2)$. We observe that $C_K(L)=L$, since it is well known that the Fitting subgroup of a solvable group is self-centralizing and $L$ is abelian. Then we have 
 \begin{equation*}
            L=C_K(L)= C_K(L_1\oplus L_2)=C_K(L_1)\cap C_K(L_2)=C_K(L_i),
        \end{equation*}
for $i=1,2$. So it follows that $L_1$ and $L_2$ are faithful and irreducible $\mathbb{F}[P]$-modules. Then, again since $P$ is cyclic, $P$ acts fixed-point free on $L_1$ and $L_2$ and therefore $P$ acts fixed-point free on $L$. So what we actually proved is that $P$ acts fixed-point free on the section $L/\Phi(L)$. We can also demonstrate that $L$ is $P\langle z \rangle$-indecomposable. Assume by contradiction that it is not true. Then there are $1\neq L_1,L_2\leq L$ such that $L_1,L_2$ are $P\langle z \rangle$-invariant (so $L_1,L_2\trianglelefteq H$) and $L=L_1\times L_2$. Then $L/\Phi(L)=L_1/\Phi(L_1)\times L_2/\Phi(L_2)$, against the fact that $L/\Phi(L)$ is minimal normal in $H/\Phi(L)$ by \Cref{prop: chief fact}. Therefore $L$ is $P\langle z \rangle$-indecomposable and $P\langle z \rangle$ acts comprimely on $L$. Then by [Corollary 1, \cite{Ha}], if $exp(L)=q^n$, we have that there exists the following normal $P\langle z \rangle$-series
 
\begin{equation*}
            L=\Omega_n(L)\trianglerighteq \Omega_{n-1}(L)\trianglerighteq \dots \trianglerighteq \Omega_0(L)=1
        \end{equation*}
such that every factor $\Omega_{n-i}(L)/\Omega_{n-i-1}(L)$ is $P\langle z \rangle$-isomorphic to $L/\Phi(L)$, for every $i=1,\dots,n-1$. Also, since $L$ is an homocyclic $q$-group, it is known that $\Omega_n(L)/\Omega_{n-1}(L)=L/\Phi(L)$. We already know that $P$ acts fixed-point free on this factor, so we conclude that the series is a normal $P$-series such that $P$ acts fixed-point free on every section. Then $P$ acts fixed-point free on $L$, that is $K=L\rtimes P$ is a Frobenius group.

We are finally ready to conclude part (a). We recall $H\cong G/N$ and so $NL\trianglelefteq G$. Then $[NL,NL]\trianglelefteq G$. Also $[NL,NL]=[N,NL][L,NL]$ and so by easy computation we get $[NL,NL]=[N,L]$. Therefore $[N,L]\trianglelefteq G$. Suppose by contradiction that $[N,L]=1$. Then every element of $N$ commutes with every element of $L$. Considering the action of $z$ on $N-\{1\}$, we see that there is a fixed-point $1\neq x\in N$ such that $x^z=x=x^{-1}$, since $N$ is elementary abelian. That is, $x$ is a real $2$-element inverted by $z$. Since $L$ contains a real non-trivial $q$-element inverted by $z$, we would get by \Cref{lem: suff cond} a real element of non prime power order, against \textbf{P}. So $1 \neq [N,L]\trianglelefteq N$ and by the minimality of $N$ follows $[N,L]=N$. Since by Fitting decomposotion we have $N=[N,L]\times C_N(L)$, then $C_N(L)=1$. So applying [Theorem 15.16, \cite{I}] to the $\mathbb{Z}_2[K]$-module $N$ we get $C_N(P)>1$. Also $C_N(P)^z=C_{N^z}(P^z)=C_N(P)$ and so $z\in N_G(C_N(P))$. Considering then the action of $z$ on $C_N(P)-\{1\}$ we see that there exists a fixed-point $1\neq x\in C_N(P)$ such that $x^z=x=x^{-1}$, that is $x$ is real inverted by $z$ and commutes with every element of $P$. But also every element of $P$ is real inverted by $z$, and so, as usual, this is against \textbf{P}. Therefore, if $3\nmid |K|$, $K$ is a $p$-group.

(b) Assume $3\mid |K|$. In this case we have that $K/\textbf{F}(H)$ is metabelian, by (2) of \Cref{theo: main DGN}. Obviously we still have $\textbf{F}(H)=\textbf{F}(K)$, $\textbf{F}(K)$ is an abelian $q$-group and $\textbf{F}(K)\in Syl_q(H)$, by \Cref{prop: Fitting is p group}. We write again $L=\textbf{F}(K)$.

Assume by contradiction that $L\nleq K'$. Then $q\mid |K:K'|$ and by \Cref{theo: structure Syl} $L$ is cyclic. But, being the Fitting subgroup of a solvable group, $L$ is self-centralizing and it is also abelian. So $C_H(L)=L$ and therefore $H/L\lesssim Aut(L) \cong Aut(C_{q^a})$, thus $H/L$ is abelian. Since $\textbf{O}^{2'}(H/L)=H/L$ we have that $H/L$ is a $2$-group, that is $H/L\cong \langle z \rangle$ and $L=K$. So $K$ would be a $q$-group, against the assumptions. Then $L\leq K'$. Since $K/L$ is metabelian, we have that $K'/L= (K/L)'$ is abelian and so $K''\leq L$. We note that we can assume $K'/L\neq 1$, otherwise we are in the case of part (a). Now we prove $K''=L$. Suppose by contradiction that $K''<L$. Then $L/K''=\textbf{F}(H/K'')$ by \Cref{prop: Fitting is p group}, since $L/K''$ is a normal Sylow $q$-subgroup of $H/K''$ and $H/K''$ satisfies (H2) and \textbf{P}. We then have, since $K'/K''$ is abelian, that $K'/K''\leq C_{H/K''}(L/K'')$, against the fact that the Fitting subgroup is self-centralizing and $L<K'$.

Now we want to verify that $K'/L$ is a $t$-group, for some (odd) prime $t$ dividing $|K|$. Let us consider $H/L$. It is easy to see that $H'=K$ and so $K'=H''$, since $H/H'$ is abelian and $\textbf{O}^{2'}(H/H')=H/H'$, that is $H/H'\cong \langle z \rangle$. Then $K'/L= (H/L)''\,char\,H/L$. Since $K'/L$ is also abelian, then $K'/L\leq \textbf{F}(H/L)$. But $\textbf{F}(H/L)$ is a $t$-group by \Cref{prop: Fitting is p group}. So $K'/L$ is a $t$-group and $t\neq q$, since $L<K'$ and $L\in Syl_q(H)$, and $t\neq p$, since $K'=K^{\infty}$ is the residual nilpotent of $K$. In conclusion, we have $K'/L\trianglelefteq H/L$ and $K'/L\in Syl_t(H/L)$. So, by \Cref{prop: Fitting is p group}, it must be $\textbf{F}(H/L)=K'/L$. 

Let $P$ be a Sylow $p$-subgroup of $H$. As in (a), we can assume $z\in N_H(P)$ and so $z$ acts fixed-point free on $P$, that is $z$ invertes every element of $P$. Also let $T$ be a Sylow $t$-subgroup of $H$ such that $TP\leq H$, that exists since $H$ has Hall $\{t,p\}$-subgroups.

\vspace{4.5cm}
\begin{figure}[ht]
\setlength{\unitlength}{0,5cm}
\begin{picture}(20,20)(-12.5,-12,5)
\put(-1,4){\line(1,0){2}}
\thicklines
\put(0,1.5){\line(0,1){13.5}}
\put(-1,15){\line(1,0){2}}
\put(1,15){\hbox{\kern3pt\ $H$}}
\put(-1,13){\line(1,0){2}}
\put(1,13){\hbox{\kern3pt\ $K=H'=H^{\infty}$}}
\put(-1,10){\line(1,0){2}}
\put(1,10){\hbox{\kern3pt\ $H''=K'=K^{\infty}$}}
\put(-1,6){\line(1,0){2}}
\put(1,6){\hbox{\kern3pt\ $H'''=K''=\textbf{F}(K)=L$}}
\put(1,4){\hbox{\kern3pt\ $\Phi(L)=\Phi(H)$}}
\put(-1,1.5){\line(1,0){2}}
\put(-3,13.75){$\langle z \rangle \cong C_2$}
\put(-3,11.25){$P\cong C_3$}
\put(-5,7.75){$T\cong C_{t^{\alpha}}\times C_{t^{\alpha}}$}
\put(-5,4.75){$C_q\times C_q \times C_q$}

\end{picture}
\vspace{-6.75cm}
\caption{Structure of group $H$ in part (b)}
\end{figure}
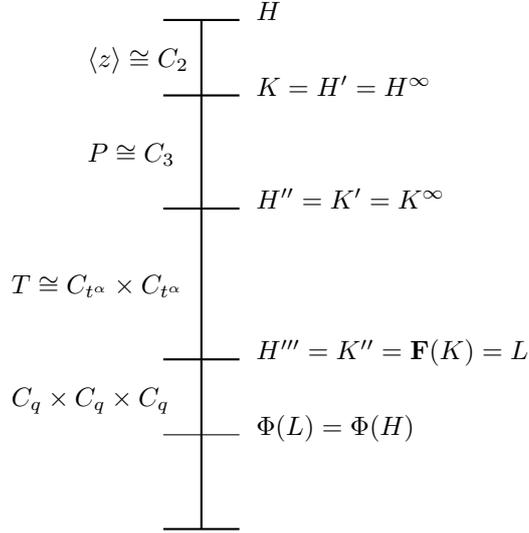 
\vspace{0.5cm}

We have that $H/L$ has the same structure and satisfies the same hypotheses as the group $H$ in part (a). Then, with the exact same reasoning, we have that $K/L$ is a Frobenius group.

Now we want to prove that $\textbf{Z}(K')=1$ so that $C_L(T)=1$, similarly to what we did in (a). Let $\Tilde{T}$ be a conjugate of $T$ such that it contains a non-trivial real element inverted by $z$, that exists by [Corollary B, \cite{DMN}]. Note that obviously $K'=L\rtimes \Tilde{T}$. Working by contradiction, assume $\textbf{Z}(K')>1$. The action of $z$ on $\textbf{Z}(K')$ is not trivial, otherwise $1<\textbf{Z}(K')\leq \textbf{Z}(H)$, against \Cref{theo: structure Syl}. So, by (4) of \Cref{lem:elem prop}, there is a real $1\neq x=ly \in \textbf{Z}(K')$, with $l\in L$, $y\in \Tilde{T}$, inverted by $z$. If $l\neq 1$, then an appropriate power of $x$ is a real $q$-element inverted by $z$ that commutes with $\Tilde{T}$, against \textbf{P}. Then $l=1$, $x\in \Tilde{T}$ and $x$ commutes with $L$, again against \textbf{P}.

We can now demonstrate that $|P|=3$. Recall $TP\cong K/L$ is a Frobenius group, with kernel $T$ and complement $P$. Consider the section $\overline{L}=L/\Phi(L)$. Since $TP$ acts coprimely on $\overline{L}$, by the previous paragraph, we have $C_{\overline{L}}(T)=\overline{C_L(T)}=1$. Also $\overline{L}$ is a $\mathbb{Z}_q$-vector space. So $\overline{L}$ is a $\mathbb{Z}_q[TP]$-module with $C_{\overline{L}}(T)=1$. Then by [Theorem 15.16, \cite{I}] we get that $|P|$ divides $\dim_{\mathbb{Z}_q}(\overline{L})=$ rank$(L)$, which is at most 3 by \Cref{theo: structure Syl}. Therefore we have $|P|=$ rank$(L)=3$.

At this point, we are able to conclude part (b). By \Cref{theo: main DGN} we have that if $N>1$ then $N=\textbf{F}(G)$. Recall also we are assuming $N$ to be elementary abelian. Therefore we have $C_G(N)=N$ and so $N$ is a faithful $G/N\cong H$-module. Let us consider the action of $TP$ on $N$. We know that $|P|=3$, $P$ acts faithfully (actually fixed-point free) on the $t$-group $T$ and $TP$ acts faithfully on $N$. In particular $N$ is a faithful $\mathbb{Z}_2[TP]$-module. Then by [36.2, \cite{Asch}] we have $C_N(P)>1$. With the same argument as in the last paragraph of (a), we get to a contradiction. Therefore, even if $3\mid |K|$, $K$ is a $p$-group.

So, in any case, $K$ is a $\langle z \rangle$-invariant $p$-section and then we have that $[K,z]$ is cyclic, since $H$ satisfies (H2). By \Cref{lem: interder KQ} $[K,z]=K$ and so $K$ is cyclic. Then, by \Cref{lem: abel K}, we have $C_K(z)=1$. Therefore $z$ invertes every element of $K$, proving that $H=K\rtimes \langle z \rangle$ is a dihedral group.

\end{proof}

\end{theorem}

Since we just did the hard work in \Cref{theo: case |Q|=2}, we can now prove Theorem A:

\begin{proof}[Proof of Theorem A]
It is obvious that (2) implies (1).

Conversely, assume that every real element of $G$ has prime power order, that is $G$ satisfies \textbf{P}. By \Cref{theo: SH} there exists $H\leq G$ such that $G=N\rtimes H$. So $H\cong G/N$. Let $Q\in Syl_2(H)$. Then, by \Cref{theo: main DGN}, we have that $H=K\rtimes Q$, where $K$ is the normal $2$-complement of $H$, and also, since $N=\textbf{O}_2(G)>1$, $Q$ is cyclic. Let $z$ be its unique involution.

Suppose $|Q|>2$. Then, still by \Cref{theo: main DGN}, we have that $K$ is cyclic. Therefore, since also $C_K(z)=C_K(Q)$, by \Cref{rem: nil impl p group} and \Cref{lem: abel K}, $K$ is a $p$-group and $C_K(z)=1$, so that every element of $K$ is real inverted by $z$. Moreover, $H=K\rtimes Q$ is a Frobenius group. Indeed, if $1\neq q \in Q$ was such that $C_K(q)>1$, since $z\in \langle q \rangle$, we would have $C_K(z)>1$, which is not true.

If $|Q|=2$, then by \Cref{theo: case |Q|=2} we have that $K$ is a cyclic $p$-group and $H$ is a dihedral group. In particular, every element of $K$ is real inverted by $z$.

In any case, we note that $G$ is a $\{2,p\}$-group and $H=K\rtimes Q$ is a Frobenius group.

Now we verify that even $N\rtimes K$ is a Frobenius group. Working by contradiction, assume that there exists $1\neq k\in K$ with $C_N(k)>1$. We see that $z$ acts on $C_N(k)$, since
\begin{equation*}
            C_N(k)^z=C_{N^z}(k^z)=C_N(k^{-1})=C_N(k).
    \end{equation*}
Then, considering the action of $z$ on $C_N(k)-\{1\}$, we see that there is a fixed-point, that is an element $1\neq x \in C_N(k)$ such that $x^z=x$. If we take an appropriate power of $x$, then we have an involution $y\in C_N(k)$ with $y^z=y=y^{-1}$. Therefore, by \Cref{lem: suff cond}, $ky$ is a real element of $G$ inverted by $z$ with non prime power order, against \textbf{P}.

So $G=NKQ$ is a $2$-Frobenius group.

Finally, let $x\in G$ be an involution. Since $K$ acts fixed-point free on $N$ and $Q$ acts fixed-point free on $K$, we have that $C_G(x)$ is a $2$-group. So $G$ is a CIT-group (see \cite{S} for the definition and more). Since $G$ is a $\{2,p\}$-group and a CIT-group, we have that every element of $G$ has prime power order.

\end{proof}

If we write $\omega_{\mathbb{R}}(G)$ to indicate the set of orders of the real elements of $G$ and if we write $\pi(\omega_{\mathbb{R}}(G))$ to denote the set of primes that divide some element of $\omega_{\mathbb{R}}(G)$, then we can state the following corollary:

\begin{corollary}\label{cor: real spectrum}
Let $G$ be a finite solvable group such that every real element has prime power order. Suppose also that $\textbf{O}_2(G)>1$. Then either $\pi(\omega_{\mathbb{R}}(G))=\{2\}$ or  $\pi(\omega_{\mathbb{R}}(G))=\{2,p\}$, for some odd prime $p$.

\begin{proof}
If $G$ is a $2$-group then $\pi(\omega_{\mathbb{R}}(G))=\{2\}$. Assume $G$ is not a $2$-group and consider $H=\textbf{O}^{2'}(G)$. We have that $H$ is a solvable group such that every real element has prime power order and $\textbf{O}^{2'}(H)=H$. Moreover, by definition, $H$ contains a Sylow $2$-subgroup of $G$ and so, since $\textbf{O}_2(G)$ is contained in every Sylow $2$-subgroup of $G$, we have $\textbf{O}_2(H)>1$. If $H$ is a $2$-group, then it is a Sylow $2$-subgroup of $G$. Since the real elements of $G$ are the real elements of $H$, then in this case  $\pi(\omega_{\mathbb{R}}(G))=\{2\}$. Lastly, suppose that $H$ is not a $2$-group. Then $H$ satisfies the hypotheses of Theorem A and so $H$ is a $\{2,p\}$-group, for some odd prime $p$. Then in this case $\pi(\omega_{\mathbb{R}}(G))=\{2,p\}$.
\end{proof}

\end{corollary}

So, given a $\{2,p\}$-group $G$ as in Theorem A, we are able to say a lot about the structure of $G/N$. But, in general, there is not much we can say about the structure of $N$. While, for a fixed $p$, we know by \cite{Hi} that there is an upper bound for dl$(N)$, the derived length of $N$, depending only on $p$, we can show that, as $p$ varies, there are groups with dl$(N)$ arbitrarily large.

\begin{example}
    This example is based on a construction of I. M. Isaacs contained in \cite{I1}. Although our results are very similar to those in \cite{I1}, since our working hypotheses are different, where necessary we will give explicit proofs.

    Let $k=2^a$, with $a\geq 1$ positive integer and let $\mathbb{F}=GF(2^k)$. Consider then $\mathcal{G}=Gal(\mathbb{F}| \mathbb{Z}_2)$, which is cyclic of order $k$, generated by $\sigma$, the Frobenius automorphism. Thus we can define $\mathbb{F}\{X\}$, the "twisted polynomial ring" in the indeterminate $X$, that is the ring of "polynomials" $\alpha_0+\alpha_1X\dots+\alpha_mX^m$ for which $X\alpha=\alpha^{\sigma}X$, for every $\alpha \in \mathbb{F}$. It is known that this does define a ring. We note that $X^k\mathbb{F}\{X\}=\mathbb{F}\{X\}X^k$ and so $(X^k)$ is a bilateral ideal. Then we consider the quotient $R=\mathbb{F}\{X\}/(X^k)$ and we write $x$ to denote the image of $X$ in $R$ under the natural homomorphism. So every element of $R$ is of the form $\alpha_0+\alpha_1x+\dots+\alpha_{k-1}x^{k-1}$ and then $|R|=(2^k)^{k}$. Also $x^k=0$ and $x\alpha=\alpha^{\sigma}x$, for every $\alpha \in \mathbb{F}$. Moreover we have that $xR=Rx$ is a nilpotent ideal and $R/xR\cong \mathbb{F}$. Therefore $xR=J(R)$, the Jacobson radical of $R$, that we denote as $J$. We have $J^i=x^iR=Rx^i$ and so $J^{k-1}\neq 0$ e $J^k=0$.  Let $S=1+J=\{1+\alpha_1x+\dots+\alpha_{k-1}x^{k-1} \mid \alpha_i \in \mathbb{F}\}$, where $1$ is the identity element of $R$. It is a known fact in ring theory that $S$ is a subgroup of the group of units of  $R$;  in this case it is a $2$-group. For $u\geq 1$ integer, we write $S_u=1+J^u$. Then $S_u\leq S$ and 
   \begin{equation*}
       S=S_1>S_2>\dots >S_{k-1}>S_k=1.
   \end{equation*}
Every element $s\in S_u$ is uniquely of the form $s=1+\alpha x^u+y$, with $y\in J^{u+1}$. We can then define $\psi_u :S_u \rightarrow \mathbb{F}$ with $\psi_u(s)=\alpha$.  It is easy to prove that $\psi_u$ is an homomorphism from $S_u$ to the additive group of $\mathbb{F}$ and that  ker$(\psi_u)=S_{u+1}$, so $S_{u+1}\trianglelefteq S_u$. By [Corollary 4.2, \cite{I1}] we have that, if $u,v\geq 1$, then $[S_u,S_v]\leq S_{u+v}$. Furthermore, if $u+v\leq k-1$, $s \in S_u$, $t \in S_v$  and $\psi_u(s)=\alpha$ and $\psi_v(t)=\beta$, then $\psi_{u+v}([s,t])=\alpha \beta^{\sigma^u}-\beta \alpha^{\sigma^v}$.

We now have to study the map $\langle \cdot,\cdot \rangle:\mathbb{F}\times \mathbb{F}\rightarrow \mathbb{F}$ such that $\langle\alpha,\beta \rangle=\alpha \beta^{\sigma^u}-\beta \alpha^{\sigma^v}$, in the case  $u$ odd. Note that $\langle \cdot,\cdot \rangle$ is $\mathbb{Z}_2$-bilinear.

We firstly prove that, if $u$ is odd and $u+v\leq k-1$, then, for $0\neq \alpha \in \mathbb{F}$,  $\langle \alpha,\mathbb{F}\rangle$ contains a hyperplane of $\mathbb{F}$. Since $\langle \alpha,\cdot\rangle$ is $\mathbb{Z}_2$-linear, it is enough to prove that its kernel  is of dimension at most $1$. So let $\langle \alpha,\beta \rangle=0=\langle \alpha,\gamma\rangle$, with $\beta\neq 0$. We then have 
   \begin{equation*}
       \alpha\beta^{\sigma^u}-\beta\alpha^{\sigma^v}=0=\alpha\gamma^{\sigma^u}-\gamma\alpha^{\sigma^v},
   \end{equation*} 
   which implies, if $\gamma\neq 0$,
   \begin{equation*}
    \beta^{-1}\beta^{\sigma^{u}}=\alpha^{-1}\alpha^{\sigma^v}=\gamma^{-1}\gamma^{\sigma^u},
   \end{equation*}
   and so $\gamma\beta^{-1}=(\gamma\beta^{-1})^{\sigma^u}$. Since $u$ is odd, we have $\langle \sigma^u\rangle=\langle \sigma \rangle$ and hence $\gamma\beta^{-1}\in \mathbb{Z}_2$.

Now we prove that there exist $0\neq\alpha,\beta \in \mathbb{F}$ such that $\langle \alpha,\mathbb{F}\rangle\neq \langle\beta,\mathbb{F}\rangle$ (even if $u$ is not odd).
To do so we will use tha trace of the Galois extension $\mathbb{F}| \mathbb{Z}_2$, that is $T:\mathbb{F}\rightarrow\mathbb{Z}_2$ with $T(\alpha)=\sum_{i=0}^{k-1} \sigma^i(\alpha)$. Recall that $T$ is invariant with respect to $\mathcal{G}=Gal(\mathbb{F}|\mathbb{Z}_2)$. Then, if $\alpha=1$, we have
   \begin{equation*}
       T(\langle 1,\gamma\rangle)=T(\gamma^{\sigma^u}-\gamma)=T(\gamma^{\sigma^u})-T(\gamma)=T(\gamma)-T(\gamma)=0.
   \end{equation*}
   Therefore it suffices  to find $\beta,\gamma\in \mathbb{F}$ with $T(\langle \beta,\gamma\rangle)\neq 0$. We have
   \begin{equation*}
       T(\langle\beta,\gamma\rangle)=T(\beta\gamma^{\sigma^u})-T(\gamma\beta^{\sigma^v})=T(\beta^{\sigma^v}\gamma^{\sigma^{u+v}})-T(\gamma\beta^{\sigma^v})=T(\beta^{\sigma^v}(\gamma^{\sigma^{u+v}}-\gamma)).
   \end{equation*}
   Since $u+v\leq k-1$ we have $k\nmid u+v$ and so we can choose $\gamma$ such that $\gamma^{\sigma^{u+v}}-\gamma\neq 0$. Moreover, simply because $\sigma^v$ is an automorphism of $\mathbb{F}$, we can choose $0\neq \beta$ so that $\beta^{\sigma^v}(\gamma^{\sigma^{u+v}}-\gamma)$ is an arbitrary element of $\mathbb{F}$, and in particular one with nonzero trace. Thus $\langle1,\mathbb{F}\rangle\neq \langle \beta,\mathbb{F}\rangle$, if $\beta$ is as said.

  In order to estimate the derived length of $S$ it is enough to verify that $[S_u,S_u]=S_{2u}$ when $u$ is odd. Indeed, assuming it is true , and denoting with $S^{(n)}$ the  $n+1$-th term of the derived series, by induction over $n$ we prove $S^{(n)}\geq S_{t_n}$, where $\{t_n\}_{n=1}^{\infty}$ is the succession defined by $t_1=2$ and $t_{n+1}=2t_n+2=2(t_n+1)$ (note that every term of $t_n$ is even and so $t_{n}+1$ is odd).
  \newline If $n=1$, we have $S^{(1)}=[S,S]=[S_1,S_1]=S_2$.
  \newline If $n>1$, then  $S^{(n+1)}=[S^{(n)},S^{(n)}]\geq [S_{t_n},S_{t_n}]\geq [S_{t_n+1},S_{t_n+1}]=S_{2t_n+2}=S_{t_{n+1}}$.
  
Therefore if we take $k$ such that $k-1\geq t_n$ we then have $S^{(n)}\geq S_{t_n}>1$ and so dl$(S)>n$.

 Let us prove $[S_u,S_u]=S_{2u}$ when $u$ is odd . Actually, fixing $u$ positive integer odd, in order to have a stronger inductive hypothesis, we prove $[S_u,S_v]=S_{u+v}$ for every $v$ positive integer. We know that $[S_u,S_v]\leq S_{u+v}$ if $u,v$ positive integers. If $u+v>k-1$, we have $S_{u+v}=1$ and we conclude. Assume then $u+v\leq k-1$. We work by induction over $(k-1-u)-v$.
  We know that $\psi_{u+v}([S_u,S_v])$ is a $\mathbb{Z}_2$-subspace of $\mathbb{F}$ which contains all the elements of the form $\langle \alpha,\beta \rangle$, with $\alpha,\beta \in \mathbb{F}$. Since $u+v\leq k-1$, we have also seen that $\psi_{u+v}([S_u,S_v])$ contains two different hyperplanes of $\mathbb{F}$ and so it has to be all of $\mathbb{F}$. Then $\psi_{u+v}([S_u,S_v])=\psi_{u+v}(S_{u+v})$.
  \newline If $(k-1-u)-v=0$, that is $u+v=k-1$, then we have ker$(\psi_{u+v})=S_{k}=1$ and so ker$(\psi_{u+v})\leq [S_u,S_v]$.
  \newline If $(k-1-u)-v>0$, we have by inductive hypothesis ker$(\psi_{u+v})=S_{u+v+1}=S_{u+(v+1)}=[S_u,S_{v+1}]\leq[S_u,S_v]$.
  \newline Therefore, since  ker$(\psi_{u+v})=S_{u+v+1}\leq [S_u,S_v]$ and $\psi_{u+v}([S_u,S_v])=\psi_{u+v}(S_{u+v})$ we conclude $[S_u,S_v]=S_{u+v}$.

In conclusion we have proved that, for a fixed positive integer $n$, if we take $k$ sufficiently large, we have dl$(S)>n$.

Consider now the semidirect product $\mathbb{F}^{\times}\rtimes \mathcal{G}$. We have $|\mathbb{F}^{\times}|=2^k-1$. Then we take a primitive divisor $p$ of $2^k-1$, that is $p\mid 2^k-1$ and $p\nmid 2^i-1$, for every $i$ positive integer, $i<k$, that exists by Zsigmondy's theorem [Theorem 6.2, \cite{MW}]. Since $p\nmid 2^{k/2}-1$, it follows $p\mid 2^{k/2}+1$. Therefore, if $P$ is a Sylow $p$-subgroup of $\mathbb{F}^{\times}$, we have $P\leq [\mathbb{F}^{\times},z]$, where $z$ is the unique involution of $\mathcal{G}$, since $\mathbb{F}^{\times}=[\mathbb{F}^{\times},z]\times C_{\mathbb{F}^{\times}}(z)$ by coprime action and $|C_{\mathbb{F}^{\times}}(z)|=2^{k/2}-1$ by Galois theory, and so $\mathcal{G}$ acts fixed-point free on $P$. Moreover, $P$ is the unique Sylow $p$-subgroup of $\mathbb{F}^{\times}\rtimes \mathcal{G}$ and so $P\trianglelefteq (\mathbb{F}^{\times}\rtimes \mathcal{G})$. We identify $P$ with the subgroup $1\cdot P$ of the unit group $R^{\times}$ of $R$. 

We now verify that $P\leq N_{R^{\times}}(S)$, that is $P$ acts on $S$, and that such action is fixed-point free. Let $s\in S$, $\gamma \in P$ and assume $s\neq 1$. So $s\in S_u$ for some $u\leq k-1$ positive integer, and we write $s=1+\alpha x^u+y$, with $0\neq \alpha \in \mathbb{F}$ and $y \in J^{u+1}$. We have 
  \begin{equation*}
      \gamma^{-1}s\gamma=1+\alpha \gamma^{-1}x^u\gamma+\gamma^{-1}y\gamma=1+\alpha \gamma^{-1}\gamma^{\sigma^u}x^u+\gamma^{-1}y\gamma.
  \end{equation*}
  It follows that $\gamma^{-1}s\gamma \in S$, that is $P$ acts on $S$, and the action of $\gamma \in P$ on $s\in S$ consist in multiplying the $u$-th coefficient of $s$ by $\gamma^{-1}\gamma^{\sigma^u}$, for $u=1,\dots,k-1$. Recall that $\sigma$ is the Frobenius automorphism of the field $\mathbb{F}$, which is of characteristic $2$, thus $\sigma(\alpha)=\alpha^2$ for every $\alpha \in \mathbb{F}$, and so $\gamma^{-1}\gamma^{\sigma^u}=\gamma^{2^u-1}$. If $1\neq s$ was fixed by $1\neq\gamma \in P$ there would be some integer $u$, $1\leq u \leq k-1$, such that the $u$-th coefficient of $s$ is nonzero and that would imply $\gamma^{2^u-1}=1$, so that $p\mid 2^u-1$, against the choice of $p$.

   So we proved that $SP$ is a Frobenius group with Frobenius kernel $S$ and complement $P$. Moreover, the automorphism $\sigma$ of  $\mathbb{F}$ can be extended to an automorphism of the ring $R$ by setting $x^{\sigma}=x$. We note that such extension, that we still denote with $\sigma$, fixes $S$ setwise and so $\mathcal{G}=\langle \sigma \rangle$ acts on $S$. Then we have the group $G=S\rtimes (P\rtimes \mathcal{G})$, that is a $\{2,p\}$-group and a $2$-Frobenius group. It easily follows then that every element of $G$ has prime power order. It is also not difficult to verify that $S=\textbf{O}_2(G)$ and $\textbf{O}^{2'}(G)=G$.
  
  In conclusion, we built a family of groups as in Theorem A, parameterized by $k=2^a$, such that, denoting $N=S=\textbf{O}_2(G)$, if $k-1\geq t_n$, we have dl$(N)>n$. So as $k$ increases we get dl$(N)$ arbitrarily large.

\end{example}

\section{Case $\textbf{O}_2(G)=1$ and consequences on $\Gamma_{\mathbb{R}}(G)$}

Let us now consider $1\neq G$ a finite solvable group, with $\textbf{O}^{2'}(G)=G$ and $\textbf{O}_2(G)=1$. Suppose that $G$ satisfies \textbf{P}. Let $Q\in Syl_2(G)$. Then by (1) of \Cref{theo: main DGN} we have that $Q$ is cyclic or generalized quaternion, $G$ has a normal $2$-complement $K$ and also $C_K(Q)=C_K(z)$, where $z$ is the unique involution of $Q$.

Regarding the structure of $\textbf{F}(G)$ we can easily deduce the following result.

\begin{lemma}\label{lem: F(G) structure}
  Let $1\neq G$ be a finite solvable group, with $\textbf{O}^{2'}(G)=G$ and $\textbf{O}_2(G)=1$. Suppose that $G$ satisfies \textbf{P}. Then there is one and only one (odd) prime $p$ such that $\textbf{O}_p(G)\nleq \textbf{Z}(G)$.

  \begin{proof}
  Firstly let us prove that such prime $p$ exists. Working by contradiction, assume it is not true. Then we have $\textbf{F}(G)\leq \textbf{Z}(G)$ and so $C_G(\textbf{F}(G))=G$. Since $G$ is solvable, $\textbf{F}(G)$ is self-centralizing and so it must be $\textbf{F}(G)=G$, which implies that $G$ is nilpotent. Since $\textbf{O}^{2'}(G)=G$, we have that $G$ is a $2$-group, against assumptions.

  Assume now by contradiction that there are two different (odd) primes $p,q$ with $\textbf{O}_p(G)\nleq \textbf{Z}(G)$ and $\textbf{O}_q(G)\nleq \textbf{Z}(G)$. Recall $Q\in Syl_2(G)$ and $z$ is its unique involution. We have that $z$ does not centralize either $\textbf{O}_p(G)$ or $\textbf{O}_q(G)$, otherwise, as already seen many times, we would get $\textbf{O}_p(G)\leq \textbf{Z}(G)$ and $\textbf{O}_q(G)\leq \textbf{Z}(G)$. Therefore there exist, by (4) of \Cref{lem:elem prop}, a non-trivial real element $x\in \textbf{O}_p(G)$ and a non-trivial real element $y\in \textbf{O}_q(G)$, both inverted by $z$. Since they commute, $xy$ is a non-trivial real element of non prime power order, against the property \textbf{P}.
  \end{proof}

\end{lemma}

It is not difficult to see that $\textbf{Z}(G)\leq K$. Also, since $G/G'$ is abelian and $\textbf{O}^{2'}(G/G')=G/G'$, we have that $G'$ has $2$-power index, that is $K\leq G'$ (actually $G'=(K\rtimes Q)'=K'[K,Q]Q'=KQ'$, since $[K,Q]=K$ by \Cref{lem: interder KQ}). So we have $\textbf{Z}(G)\leq G'$ and therefore $\textbf{Z}(G)\leq \Phi(G)$, but this by itself obviously does not force $\textbf{Z}(G)=1$. If it were $\textbf{Z}(G)=1$ then we could infer that $\textbf{F}(G)$ is a $p$-group.

\medskip

As already announced in the introduction, in this paper we were not able to precisely describe the structure of these groups, but we concentrated our efforts on finding the following examples, that prove that in this case $G$ is not necessarily a $\{2,p\}$-group.

\begin{example}\label{ex: G non 2 p group}
Consider $S_3$ and its unique, up to isomorphism, irreducible $\mathbb{Z}_5[S_3]$-module $V$ of dimension $2$. Then we can consider the semidirect product $G=V\rtimes S_3$. In \texttt{GAP} $G$ is the \texttt{SmallGroup(150,5)}. Obviously $G$ is solvable and its structure is $(C_5\times C_5)\rtimes S_3$. We can verify that:
\begin{itemize}
    \item the normal subgroups of $G$ are $1$, $G''=\textbf{F}(G)\cong (C_5\times C_5)$, $G'\cong (C_5\times C_5)\rtimes C_3$, and $G$. So we have $\textbf{O}_2(G)=1$ and $\textbf{O}^{2'}(G)=G$;
    \item the orders of the real elements of $G$ are $\{1,2,3,5\}$, that is $|\pi(\omega_{\mathbb{R}}(G))|=3$ and every real element has prime (power) order;
    \item $\textbf{Z}(G)=1$ (consistenly with the fact that $\textbf{F}(G)$ is a $5$-group).
\end{itemize}

Note that the structure of $G$ resembles the one in Figure 1.

\medskip

We can also extend $G$ by considering the semidirect product with an appropriate irreducible $\mathbb{Z}_{11}[G]$-module $W$ of dimension $3$. Unfortunately, since the resulting group $H=W\rtimes G$ has order 199650, it is not listed in \texttt{GAP} libraries. So we are gonna identify $W$, and so also $H$, by saying that $G$ acts on $W$, as a group of matrices, with generators: 
\begin{equation*}
     \begin{bmatrix}
         \bar{0} & \bar{0} & \bar{1}\\
         \bar{0} & \bar{1} & \bar{0}\\
         \bar{1} & \bar{0} & \bar{0}
     \end{bmatrix}
     \begin{bmatrix}
        \bar{0} & \bar{0} & \bar{1}\\
         \bar{1} & \bar{0} & \bar{0}\\
         \bar{0} & \bar{1} & \bar{0}
     \end{bmatrix}
     \begin{bmatrix}
          \bar{9} & \bar{0} & \bar{0}\\
         \bar{0} & \bar{3} & \bar{0}\\
         \bar{0} & \bar{0} & \bar{9}
     \end{bmatrix}
     \begin{bmatrix}
          \bar{5} & \bar{0} & \bar{0}\\
         \bar{0} & \bar{9} & \bar{0}\\
         \bar{0} & \bar{0} & \bar{1}
     \end{bmatrix},
 \end{equation*}
 respectively of order $2,3,5,5$. Obviously $H$ is solvable and its structure is $(C_{11}\times C_{11}\times C_{11})\rtimes ((C_5\times C_5)\rtimes S_3)$. We can verify that:
 \begin{itemize}
     \item the normal subgroups of $H$ are $1$, $H'''=\textbf{F}(H)\cong C_{11}\times C_{11}\times C_{11}$, $H''= (C_{11}\times C_{11}\times C_{11})\rtimes (C_5\times C_5)$, $H'=(C_{11}\times C_{11}\times C_{11})\rtimes ((C_5\times C_5)\rtimes C_3)$ e $H$. So we have $\textbf{O}_2(H)=1$ and $\textbf{O}^{2'}(H)=H$;
     \item the orders of the real elements of $H$ are $\{1,2,3,5,11\}$, that is $|\pi(\omega_{\mathbb{R}}(G))|=4$ and every real element has prime (power) order;
     \item $\textbf{Z}(H)=1$ (consistenly with the fact that $\textbf{F}(G)$ is a $11$-group).
 \end{itemize}

 Note that the structure of $H$ resembles the one in Figure 2.

\end{example}

It could be an interesting topic for future research to see if there are further extensions and how much it is possible to increase $|\pi(\omega_{\mathbb{R}}(G))|$, preserving the hypothesis that every real element has prime power order.

\bigskip

We will now, for our final considerations, switch to the framework of the real prime graph $\Gamma_{\mathbb{R}}(G)$. It is clear that the existence of groups $G$ and $H$ as in \Cref{ex: G non 2 p group} proves Theorem B.

Let us also consider other known graphs associated to $G$ and then compare their properties with the ones of $\Gamma_{\mathbb{R}}(G)$.

Regarding the Gruenberg-Kegel graph $\Gamma(G)$, we have to recall the following simple, but still extremely important, result by M. S. Lucido, also known as Lucido's Lemma (see [Proposition 1, \cite{L}]):

\begin{lemma}\label{lem: Lucido}
Let $G$ be a finite solvable group. If $p,q,r$ are three distinct primes that divide $|G|$, then $G$ contains an element of order the product of two of these primes.
    
\end{lemma}

It is immediate that \Cref{lem: Lucido} is equivalent to the fact that $\Gamma(G)$ does not contain a set of three pairwise non-adjacent vertices, and so in particular $n(\Gamma(G))\leq 2$. Since the group $H$ (also $G$) from the \Cref{ex: G non 2 p group} is such that $\Gamma_{\mathbb{R}}(H)$ contains sets of three pairwise non-adjacent vertices, we deduce that an analogous version of \Cref{lem: Lucido} for real elements does not hold.

Beyond the comparison with $\Gamma(G)$, as we already mentioned, we can compare this feature of $\Gamma_{\mathbb{R}}(G)$ with other two notorius graphs associated to $G$, that historically have maintained certain similarities with respect to $\Gamma_{\mathbb{R}}(G)$. Such graphs are $\Gamma_{cd,\mathbb{R}}(G)$, the prime graph on real character degrees, and $\Gamma_{cs,\mathbb{R}}(G)$, the prime graph on real class sizes.

Regarding the least upper bound for the number of connected components of these graphs, we have the following results, respectively [(ii) of Theorem 5.1, \cite{DNT}] and [Theorem 6.2, \cite{DNT}]:

\begin{theorem}
    Let $G$ be a finite solvable group. Then $n(\Gamma_{cd,\mathbb{R}}(G))\leq 2$.
\end{theorem}

\begin{theorem}
    Let $G$ be a finite group. Then $n(\Gamma_{cs,\mathbb{R}}(G))\leq 2$.
\end{theorem}

So we deduce that, at least in this case, there is a breaking of the symmetry between $\Gamma_{cs,\mathbb{R}}(G)$, $\Gamma_{cd,\mathbb{R}}(G)$ and $\Gamma_{\mathbb{R}}(G)$. Furthermore, we note that the least upper bound for $n(\Gamma_{\mathbb{R}}(G))$, if it exist, has to be at least $4$. So determining such bound, even though now there is a little bit more information about it, remains an open question of interest.

\section*{Acknowledgements}
The author would like to thank Professor S. Dolfi for bringing this problem to his attention and for the invaluable advice received while working under his supervision on this research.

\nocite{*}
\printbibliography

@article{DGN,
    author = {S. Dolfi and D. Gluck and G. Navarro} ,
    title = {On the orders of real elements of solvable groups},
    journal = {Israel Journal of Mathematics},
    volume= {210},
    pages= {1--21},
    year = {2015},
    %shorthand= {DGN},
}

@article{DMN,
    author = {S. Dolfi and G. Malle and G. Navarro},
    title = {The finite groups with no real p-elements},
    journal = {Israel Journal of Mathematics},
    volume= {192},
    pages= {831--840},
    year = {2012},
    %shorthand={DMN},
}

@book{I,
    author ={I. M. Isaacs} ,
    title ={Character Theory of Finite Groups} ,
    publisher ={Academic press, New York} ,
    year = {1976},
    %shorthand={I}
}

@book{Asch,
    author = {M. Aschbacher} ,
    title = {Finite Group Theory},
    publisher = {Cambridge University Press},
    edition={Second edition},
    year = {2000},
    %shorthand={Asch},
}

@article{Ha,
    author ={M. Harris} ,
    title ={On $p'$-automorphisms of abelian p-groups},
    journal = {Rocky Mountain Journal of Mathematics},
    volume= {7},
    pages={751--752},
    year = {1977},
    %shorthand= {Ha}
}

@article{S,
    author = {M. Suzuki} ,
    title = {Finite groups with nilpotent centralizers},
    journal = {Transactions of the American Mathematical Society},
    volume={99},
    pages={425--470},
    year ={1961},
    %shorthand={S},
}

@article{L,
    author = {M. S. Lucido},
    title = {The diameter of the prime graph of a finite group},
    journal = {Journal of Group Theory},
    volume={2},
    pages={157--172},
    year = {1999},
    %shorthand={L}
}

@article{DNT,
    author ={S. Dolfi and G. Navarro and P. H. Tiep} ,
    title = {Primes dividing the degrees of the real characters} ,
    journal = {Mathematische Zeitschrift},
    volume={259},
    pages= {755--774},
    year = {2008},
    %shorthand= {DNT},
}

@manual{GAP4,
    organization= {The GAP Group},
    title = {Groups, Algorithms, and Programming, Version 4.12.1},
    year={2022},
    url={https://www.gap-system.org},
    %shorthand={GAP},
}

@article{Hi,
    author ={G. Higman},
    title ={Finite groups in which every element has prime power order},
    journal={Journal of the London Mathematical Society},
    volume={S1-32},
    pages={335--342},
    year={1957},
    %shorthand={Hi},
}

@article{I1,
  title={Coprime group actions fixing all nonlinear irreducible characters},
  author={I. M. Isaacs},
  journal={Canadian Journal of Mathematics},
  volume={41},
  pages={68--82},
  year={1989},
  %shorthand={I1},
}

@book{MW,
    author={O. Manz and T. R. Wolf},
    title={Representations of Solvable Groups},
    edition={London Mathematical Society Lecture Note Series, Vol. 185},
    publisher={Cambridge University Press, Cambridge},
    year={1993},
    %shorthand={MW},
}

\end{document}